\documentclass[11pt,reqno, twoside]{article}
\usepackage{amsfonts}
\usepackage{amsfonts, amsmath, amssymb, amscd, amsthm}
\def\qedbox{\hbox{$\rlap{$\sqcap$}\sqcup$}}
\def\qed{\nobreak\hfill\penalty250 \hbox{}\nobreak\hfill\qedbox}

\usepackage{mathrsfs}
\usepackage{multicol}
\usepackage{comment}

 \newtheorem{thm}{Theorem}[section]
 
 \newtheorem{lem}[thm]{Lemma}
 
 \theoremstyle{definition}
 
 \newtheorem{ex}[thm]{Example}
 \newtheorem{rem}[thm]{Remark}
 \theoremstyle{remark}

 \theoremstyle{claim}
 \numberwithin{equation}{section}

\numberwithin{equation}{section}

\newcounter{rom}
\renewcommand{\therom}{(\roman{rom})}
{\end{list}}

\title{Second eigenvalue of a Jacobi
operator of hypersurfaces with constant scalar curvature}

\author{Haizhong Li \thanks{Supported by Tsinghua
University--K.U.Leuven Bilateral Scientific Cooperation Fund.}
\thanks {Supported by NSFC grant No. 10971110.}
\and Xianfeng Wang \protect\footnotemark[1] \thanks {Supported by
NSFC grant No. 10701007.} }

\date{}

\begin{document}

\maketitle

%%% ----------------------------------------------------------------------
\begin{abstract}
Let $x:M\to\mathbb{S}^{n+1}(1)$ be an n-dimensional compact
hypersurface with constant scalar curvature $n(n-1)r,~r\geq 1$, in a
unit sphere $\mathbb{S}^{n+1}(1),~n\geq 5$. We know that such
hypersurfaces can be characterized as critical points for a
variational problem of the integral $\int_MH dv$ of the mean
curvature $H$. In this paper, we derive an optimal upper bound for
the second eigenvalue of the Jacobi operator $J_s$ of $M$. Moreover,
when $r>1$, the bound is attained if and only if $M$ is totally
umbilical and non-totally geodesic, when $r=1$, the bound is
attained if $M$ is the Riemannian product
$\mathbb{S}^{m}(c)\times\mathbb{S}^{n-m}(\sqrt{1-c^2}),~1\leq m\leq
n-2,~c=\sqrt{\frac{(n-1)m+\sqrt{(n-1)m(n-m)}}{n(n-1)}}$.
\end{abstract}



%%% ----------------------------------------------------------------------
\maketitle
%%% ----------------------------------------------------------------------
%\tableofcontents
\section{Introduction}

Let $M$ be an n-dimensional compact hypersurface in a unit sphere
$\mathbb{S}^{n+1}(1)$. We denote the components of the second
fundamental form of $M$ by $h_{ij}$, and denote the principal
curvatures of $M$ by $k_1,\ldots,k_n$. Let $H$, $H_2$ and $H_3$
denote the mean curvature, the 2nd mean curvature and the 3rd mean
curvature of $M$ respectively, namely,
$$
\begin{aligned}
H&=\frac{1}{n}\sum_{i=1}^{n}k_i,~H_2=\frac{2}{n(n-1)}\sum_{1\leq i_1<i_2\leq n}k_{i_1}k_{i_2},\\
H_3&=\frac{6}{n(n-1)(n-2)}\sum_{1\leq i_1<i_2<i_3\leq
n}k_{i_1}k_{i_2}k_{i_3}.\end{aligned}
$$

We denote the square norm of the second fundamental form of $M$ by
$S$. The Schr\"{o}dinger operator $J_m=-\Delta-S-n,$ where $\Delta$
stands for the Laplace-Beltrami operator, is called the Jacobi
operator. Its spectral behavior is directly related to the
instability of both the minimal hypersurfaces and the hypersurfaces
with constant mean curvature in $\mathbb{S}^{n+1}(1)$ (cf. \cite{S}
and \cite{BCE}). The first eigenvalue of the Jacobi operator $J_m$
of such hypersurfaces in $\mathbb{S}^{n+1}(1)$ was studied by Simons
\cite{S} and Wu \cite{W}.

The second eigenvalue of the Jacobi operator $J_m$ of the compact
hypersurfaces in $\mathbb{S}^{n+1}(1)$ was studied by A. El Soufi
and S. Ilias in \cite{SI}. They obtained that if $M$ is an
n-dimensional compact hypersurface in $\mathbb{S}^{n+1}(1)$, then
the second eigenvalue $\lambda_2^{J_m}$ of the Jacobi operator $J_m$
satisfies
$$\lambda_2^{J_m}\leq 0,$$
where the equality holds if and only if $M$ is a totally umbilical
hypersurface in $\mathbb{S}^{n+1}(1)$.

For any $C^2$-function $f$ on $M$, we define a differential operator
\begin{equation}\Box f=\sum_{i,j=1}^{n}(nH\delta_{ij}-h_{ij})f_{ij},\label{1-1}\end{equation}
where $(f_{ij})$ is the Hessian of $f$. The differential operator
$\Box$ is self-adjoint and it was introduced by S. Y. Cheng and Yau
in \cite{CY} in order to study the compact hypersurfaces with
constant scalar curvature in $\mathbb{S}^{n+1}(1)$. They proved that
if $M$ is an n-dimensional compact hypersurface with constant scalar
curvature $n(n-1)r,~r\geq 1$, and if the sectional curvature of $M$
is non-negative, then $M$ is either a totally umbilical hypersurface
$\mathbb{S}^{n}(c)$ or a Riemannian product
$\mathbb{S}^{m}(c)\times\mathbb{S}^{n-m}(\sqrt{1-c^2}),~1\leq m\leq
n-1$, where $\mathbb{S}^{k}(c)$ denotes a sphere of radius $c$. In
\cite{L}, the first author proved that if $M$ is an n-dimensional
$(n\geq3)$ compact hypersurface with constant scalar curvature
$n(n-1)r,~r\geq 1$, and if $S\leq
(n-1)\frac{n(r-1)+2}{n-2}+\frac{n-2}{n(r-1)+2}$, then $M$ is either
a totally umbilical hypersurface or a Riemannian product
$\mathbb{S}^{1}(c)\times\mathbb{S}^{n-1}(\sqrt{1-c^2})$ with
$0<1-c^2=\frac{n-2}{nr}\leq \frac{n-2}{n}$. Furthermore, the
Riemannian product
$\mathbb{S}^{1}(c)\times\mathbb{S}^{n-1}(\sqrt{1-c^2})$ has been
characterized in \cite{Cheng2} and \cite{Cheng3}.

In \cite{ACC}, Alencar, do Carmo and Colares studied the stability
of the hypersurfaces with constant scalar curvature in
$\mathbb{S}^{n+1}(1)$. In this case, the Jacobi operator $J_s$ is
given by (cf. \cite{ACC} and \cite{Cheng4})
\begin{equation}
J_s=-\Box-\{n(n-1)H+nHS-f_3\},\label{1.1}
\end{equation}
which is associated with the variational characterization of the
hypersurfaces with constant scalar curvature in
$\mathbb{S}^{n+1}(1)$, where $f_3=\sum\limits_{j=1}^{n}k_j^3$ (cf.
\cite{Re} and \cite{Ro}). The spectral behavior of $J_s$ is directly
related to the instability of the hypersurfaces with constant scalar
curvature.

In general, $J_s$ is not an elliptic operator. When $r>1$,
$n^2H^2>S>0$, the differential operator $\Box$ and hence $J_s$ is an
elliptic operator (cf. pages 3310, 3311 in \cite{Cheng4}). When
$r=1$, if we assume that $H_3\neq 0$ on $M$, then we have $H\neq 0$
and $J_s$ is elliptic (cf. Proposition 1.5 in \cite{HL}).

\noindent{\bf Definition 1:}  We call $\lambda_i^{J_s}$ an
eigenvalue of $J_s$ if there exists a non-zero function $f$ on $M$
such that $J_sf=\lambda_i^{J_s} f$, we call $\lambda_i^{\Box}$ an
eigenvalue of $\Box$ if there exists a non-zero function $f$ on $M$
such that $\Box f+\lambda_i^{\Box} f=0,$  and we call
$\lambda_i^{\Delta}$ an eigenvalue of $\Delta$ if there exists a
non-zero function $f$ on $M$ such that $\Delta f+\lambda_i^{\Delta}
f=0.$

In \cite{Cheng4}, Q. -M. Cheng studied the first eigenvalue of $J_s$
of the hypersurfaces  with constant scalar curvature $n(n-1)r,r>1$
in $\mathbb{S}^{n+1}(1)$, and derived an optimal upper bound for the
first eigenvalue of $J_s$.
\begin{thm} (see Corollary 1.2 in \cite{Cheng4})
Let $M$ be an n-dimensional compact orientable hypersurface with constant scalar
curvature $n(n-1)r,~r>1$, in $\mathbb{S}^{n+1}(1)$. Then the Jacobi operator $J_s$ is
elliptic and the first eigenvalue of $J_s$ satisfies
$$\lambda_1^{J_s}\leq -n(n-1)r\sqrt{r-1},$$ where the equality holds if
and only if $M$ is totally umbilical and non-totally geodesic.
\end{thm}

In \cite{ABS}, L. J. Al\'{\i}as, A. Brasil and L. A. M. Sousa
studied the first eigenvalue $\lambda_1^{J_s}$ of $J_s$ of the
hypersurfaces with constant scalar curvature $n(n-1)$ in
$\mathbb{S}^{n+1}(1)$.
\begin{thm} (see Theorem 2 in \cite{ABS})
Let $M$ be an n-dimensional compact orientable hypersurface with
constant scalar curvature $n(n-1)$, in
$\mathbb{S}^{n+1}(1),~n\geq3$. Assume that $H_3\neq 0$, then the
Jacobi operator $J_s$ is elliptic and the first eigenvalue
$\lambda_1^{J_s}$ of the Jacobi operator $J_s$ satisfies
$$\lambda_1^{J_s}\leq -2n(n-1)\min{|H|},$$ where the equality holds if
and only if $M$ is the Riemannian product
$\mathbb{S}^{m}(c)\times\mathbb{S}^{n-m}(\sqrt{1-c^2})$ with $1\leq
m\leq n-2,~c=\sqrt{\frac{(n-1)m+\sqrt{(n-1)m(n-m)}}{n(n-1)}}$.
\end{thm}

In this paper, we study the second eigenvalue for $J_s$ of the
hypersurfaces  with constant scalar curvature $n(n-1)r,r\geq1$ in
$\mathbb{S}^{n+1}(1),~n\geq 5$, and we have the following results.

\begin{thm}
Let $M$ be an n-dimensional compact orientable hypersurface with
constant scalar curvature $n(n-1)r,~r> 1$, in
$\mathbb{S}^{n+1}(1),~n\geq5$. Then, the Jacobi operator $J_s$ is
elliptic and the second eigenvalue $\lambda_2^{J_s}$ of the Jacobi
operator $J_s$ satisfies $$\lambda_2^{J_s}\leq 0,$$ where the
equality holds if and only if  $M$ is totally umbilical and
non-totally geodesic.
\end{thm}

\begin{thm}
Let $M$ be an n-dimensional compact orientable hypersurface with
constant scalar curvature $n(n-1)$, in
$\mathbb{S}^{n+1}(1),~n\geq5$. Assume that $H_3\neq 0$, then the
Jacobi operator $J_s$ is elliptic and the second eigenvalue
$\lambda_2^{J_s}$ of the Jacobi operator $J_s$ satisfies
\begin{equation}\lambda_2^{J_s}\leq -\frac{
n(n-1)(n-2)}{2}\min{|H_3|},\label{1.2}\end{equation} where the
equality holds if and only if $H_3=\text{constant}\neq 0$ and the
position functions of $M$ in $\mathbb{S}^{n+1}(1)$ are the second
eigenfunctions of $J_s$ corresponding to $\lambda_2^{J_s}$. In
particular, when $M$ is the Riemannian product
$\mathbb{S}^{m}(c)\times\mathbb{S}^{n-m}(\sqrt{1-c^2}),~1\leq m\leq
n-2,~c=\sqrt{\frac{(n-1)m+\sqrt{(n-1)m(n-m)}}{n(n-1)}}$, the
equality in \eqref{1.2} is attained.
\end{thm}
\section{Preliminaries}
Throughout this paper, all manifolds are assumed to be smooth and
connected without boundary.  Let $x:M\to \mathbb{S}^{n+1}(1)$ be an
n-dimensional hypersurface in a unit sphere $\mathbb{S}^{n+1}(1)$.
We make the following convention on the range of indices:
$$1\leq i,j,k,l\leq n.$$

Let $\{e_1,\cdots,e_n,e_{n+1}\}$ be a local orthonormal frame with
dual coframe $\{\omega_1,\cdots,\omega_n,\omega_{n+1}\}$ such that
when restricted on $M$,  $\{e_1,\cdots,e_n\}$ is a local orthonormal
frame on $M$. Hence we have $\omega_{n+1}=0$ on $M$ and we have the
following structure equations (see \cite{CL}, \cite{C}, \cite{L} and
\cite{S}):
\begin{equation}
dx=\sum_{i}\omega_ie_i,\label{2.1}
\end{equation}
\begin{equation}
de_i=\sum_{j}\omega_{ij}e_j+\sum_{j}h_{ij}\omega_{j}e_{n+1}-\omega_ix,\label{2.2}
\end{equation}
\begin{equation}
de_{n+1}=-\sum_{i,j}h_{ij}\omega_je_i,\label{2.3}
\end{equation}
where $h_{ij}$ denote the components of the second fundamental form
of $M$.

The Gauss equations are (see \cite{C}, \cite{L})
\begin{equation}
R_{ijkl}=\delta_{ik}\delta_{jl}-\delta_{il}\delta_{jk}+h_{ik}h_{jl}-h_{il}h_{jk},\label{2.4}
\end{equation}
\begin{equation}
R_{ik}=(n-1)\delta_{ik}+nHh_{ik}-\sum_{j}h_{ij}h_{jk},\label{2.5}
\end{equation}
\begin{equation}
R=n(n-1)r=n(n-1)+n^2H^2-S,\label{2.6}
\end{equation}
where $R$ is the scalar curvature of $M$, $r$ is the normalized scalar curvature of $M$ and $S=\sum\limits_{i,j}h_{ij}^2$
is the norm square of the second fundamental form,
$H=\frac{1}{n} \sum\limits_{i}h_{ii}$ is the mean curvature of $M$.

The Codazzi equations are given by (see \cite{C}, \cite{L})
\begin{equation}
h_{ijk}=h_{ikj}.\label{2.7}
\end{equation}

Let $f$ be a smooth function on $M$, we define its gradient and
Hessian by (see \cite{C}, \cite{L})
\begin{equation}
df=\sum_{i=1}^{n}f_i\omega_i,\label{2.8}
\end{equation}
\begin{equation}
\sum_{j=1}^{n}f_{ij}\omega_j=df_i+\sum_{j=1}^{n}f_j\omega_{ji}.\label{2.9}
\end{equation}

Then the Jacobi operator $J_s$ (see \eqref{1.1}) is defined by
\begin{equation}
\begin{aligned}
J_sf&=-\Box f-\{n(n-1)H+nHS-f_3\}f\\
&=-\sum_{i,j}(nH\delta_{ij}-h_{ij})f_{ij}-\{n(n-1)H+nHS-f_3\}f.
\end{aligned}\label{2.10}
\end{equation}

\section{Some examples and some lemmas}
First of all, we consider the first and second eigenvalues of the Jacobi operator $J_s$
of the totally umbilical and non-totally geodesic hypersurface in $\mathbb{S}^{n+1}(1)$
with constant scalar curvature $n(n-1)r,~r>1$ and the Riemannian product
$\mathbb{S}^{m}(c)\times\mathbb{S}^{n-m}(\sqrt{1-c^2}),~1\leq m\leq n-2$ with constant
scalar curvature $n(n-1)$ in $\mathbb{S}^{n+1}(1),~n\geq 3$.

\begin{ex}
Let $M$ be a totally umbilical and non-totally geodesic hypersurface
with constant scalar curvature $n(n-1)r,~r>1$ in
$\mathbb{S}^{n+1}(1)$. We can assume $H>0$. In this case,
$\Box=(n-1)H\Delta$, from $S=nH^2$ and the Gauss equation
\eqref{2.6} we have $H=\sqrt{r-1}$. By \eqref{1.1} we have
$$J_s=-\Box-\{n(n-1)H+nHS-f_3\}=-\{(n-1)H\Delta+n(n-1)H(1+H^2)\},$$
hence the eigenvalues $\lambda_i^{J_s}$ of $J_s$ are given by
$$\lambda_i^{J_s}=(n-1)H\lambda_i^{\Delta}-n(n-1)H(1+H^2),$$
where $\lambda_i^{\Delta}$ denotes the eigenvalue of $\Delta$ (see
Definition 1). It is well-known that
$\lambda_1^{\Delta}=0,~\lambda_2^{\Delta}=nr=n(1+H^2)$, hence we
have
\begin{equation}
\begin{aligned}
\lambda_1^{J_s}&=-n(n-1)H(1+H^2)=-n(n-1)r\sqrt{r-1}<0,\\
\lambda_2^{J_s}&=(n-1)H\cdot n(1+H^2)-n(n-1)H(1+H^2)=0.
\end{aligned}
\end{equation}
\end{ex}

\begin{ex}
Let $M$ be the Riemannian product
$$\mathbb{S}^{m}(c)\times\mathbb{S}^{n-m}(\sqrt{1-c^2}),~1\leq m\leq
n-2,~c=\sqrt{\frac{(n-1)m+\sqrt{(n-1)m(n-m)}}{n(n-1)}}$$  in $\mathbb{S}^{n+1}(1),~n\geq 3$.
In this case, the position vector is
$$x=(x_1,x_2)\in
\mathbb{S}^{m}(c)\times\mathbb{S}^{n-m}(\sqrt{1-c^2})$$ and the unit
normal vector at this point $x$ is given by
$e_{n+1}=(\frac{\sqrt{1-c^2}}{c}x_1,-\frac{c}{\sqrt{1-c^2}}x_2)$.

Its principal curvatures are given by
\begin{equation}
k_1=\cdots=k_{m}=-\frac{\sqrt{1-c^2}}{c},~k_{m+1}=\cdots=k_n=\frac{c}{\sqrt{1-c^2}}.
\end{equation}

Since the principal curvatures are constant hence $H,~S,~f_3$ are
all constant given by
\begin{equation}
\begin{aligned}
H&=\frac{nc^2-m}{cn\sqrt{1-c^2}},\\
S&=\frac{m(1-c^2)}{c^2}+\frac{(n-m)c^2}{1-c^2}=n^2H^2,\\
f_3&=-\frac{m(1-c^2)^{3/2}}{c^3}+\frac{(n-m)c^3}{(1-c^2)^{3/2}}.
\end{aligned}
\end{equation}

After a long but straightforward computation, we know that $M$ has constant scalar
curvature $n(n-1)$ and
\begin{equation}
\begin{aligned}
H_3=-\frac{2H}{n-2}=-\frac{2(nc^2-m)}{cn(n-2)\sqrt{1-c^2}}<0,
\end{aligned}
\end{equation}
hence the Jacobi operator $J_s$ is elliptic (cf. Proposition 1.5 in \cite{HL}). We also
have
\begin{equation}
n(n-1)H+nHS-f_3=\frac{(n-2m)(n-1)c^4+2m(m-1)c^2-m(m-1)}{c^3(1-c^2)^{3/2}}, \label{0037}
\end{equation}
thus the Jacobi operator $J_s=-\Box-\{n(n-1)H+nHS-f_3\}$ becomes
\begin{equation}
J_s=-\Box-\frac{(n-2m)(n-1)c^4+2m(m-1)c^2-m(m-1)}{c^3(1-c^2)^{3/2}}, \label{0038}
\end{equation}
hence, the eigenvalues $\lambda_i^{J_s}$ of $J_s$ are given by
\begin{equation}
\lambda_i^{J_s}=\lambda_i^{\Box}-\frac{(n-2m)(n-1)c^4+2m(m-1)c^2-m(m-1)}{c^3(1-c^2)^{3/2}},\label{0039}
\end{equation}
where $\lambda_i^{\Box}$ denotes the eigenvalue of the differential
operator $\Box$ (see Definition 1).

Since the differential operator $\Box$ is self-adjoint and $M$ is compact, we have
$\lambda_1^{\Box}=0$ and its corresponding eigenfunctions are non-zero constant
functions, hence
\begin{equation}
\lambda_1^{J_s}=-\frac{(n-2m)(n-1)c^4+2m(m-1)c^2-m(m-1)}{c^3(1-c^2)^{3/2}}.\label{00310}
\end{equation}

Let $\{e_1,\cdots,e_n\}$ be a local orthonormal basis of $TM$ with
dual basis $\{\omega_1,\cdots,\omega_n\}$ such that
$\{e_1,\cdots,e_m\}$ is a local orthonormal basis of
$T\mathbb{S}^{m}(c)$ when restricted on $\mathbb{S}^{m}(c)$ and
$\{e_{m+1},\cdots,e_n\}$ is a local orthonormal basis of
$T\mathbb{S}^{n-m}(\sqrt{1-c^2})$ when restricted on
$\mathbb{S}^{n-m}(\sqrt{1-c^2})$. So we have
\begin{equation}\Box f=\sum_{i=1}^{m}(nH-k_1)f_{ii}+\sum_{j=m+1}^{n}(nH-k_n)f_{jj}=
(nH-k_1)\Delta_1f+(nH-k_n)\Delta_2f,\label{00311}\end{equation}
where $\Delta_1$ and $\Delta_2$ denote the Laplace-Beltrami
operators on $\mathbb{S}^{m}(c)$ and
$\mathbb{S}^{n-m}(\sqrt{1-c^2})$ respectively. Since
$(nH-k_1)=\frac{(n-1)c^2-(m-1)}{c\sqrt{1-c^2}}>0,~(nH-k_n)=\frac{(n-1)c^2-m}{c\sqrt{1-c^2}}>0$,
we conclude that
\begin{equation}
\lambda_2^{\Box}=\min{\{(nH-k_1)\lambda_2^{\Delta_1},(nH-k_n)\lambda_2^{\Delta_2}\}},\label{00312}
\end{equation}
where $\lambda_2^{\Delta_1}$ and $\lambda_2^{\Delta_2}$ are the second eigenvalues (or
the first non-zero eigenvalue) of $\Delta_1$ and $\Delta_2$ which are given by
\begin{equation}
\lambda_2^{\Delta_1}=\frac{m}{c^2},~\lambda_2^{\Delta_2}=\frac{n-m}{1-c^2}.\label{00313}
\end{equation}
Therefore, from  \eqref{00312} and \eqref{00313}, after a direct
computation, we have
\begin{equation}
\begin{aligned}
\lambda_2^{J_s}&=\text{min}\{(nH-k_1)\frac{m}{c^2}-\frac{(n-2m)(n-1)c^4+2m(m-1)c^2-m(m-1)}{c^3(1-c^2)^{3/2}},\\
&(nH-k_n)\frac{n-m}{1-c^2}-\frac{(n-2m)(n-1)c^4+2m(m-1)c^2-m(m-1)}{c^3(1-c^2)^{3/2}}\}\\
&=\min{\{\frac{(n-m)[(1-n)c^2+m]}{c(1-c^2)^{3/2}},\frac{-m[(n-1)c^2-(m-1)]}{c^3(1-c^2)^{1/2}}\}}.\label{0.314}
\end{aligned}
\end{equation}

Since $c=\sqrt{\frac{(n-1)m+\sqrt{(n-1)m(n-m)}}{n(n-1)}}$, we have
\begin{equation}
\begin{aligned}
&\frac{(n-m)[(1-n)c^2+m]}{c(1-c^2)^{3/2}}-\frac{m[(n-1)c^2-(m-1)]}{c^3(1-c^2)^{1/2}}\\
&=-\frac{n(n-1)c^4+2m(1-n)c^2+m(m-1)}{c^3(1-c^2)^{3/2}}=0.
\end{aligned}\label{312}
\end{equation}
It follows from \eqref{0.314} and \eqref{312} that
\begin{equation}
\begin{aligned}
\lambda_2^{J_s}=\frac{(n-m)[(1-n)c^2+m]}{c(1-c^2)^{3/2}}<0.
\end{aligned}\label{3-30}
\end{equation}
On the other hand, we also have
\begin{equation}
\begin{aligned}
&-\frac{(n-2m)(n-1)c^4+2m(m-1)c^2-m(m-1)}{c^3(1-c^2)^{3/2}}+2n(n-1)H\\
&=-\frac{(2c^2-1)(n(n-1)c^4+2m(1-n)c^2+m(m-1))}{c^3(1-c^2)^{3/2}}=0,
\end{aligned}\label{316}
\end{equation}
\begin{equation}
\begin{aligned}
&\frac{(n-m)[(1-n)c^2+m]}{c(1-c^2)^{3/2}}+n(n-1)H\\
&=-\frac{n(n-1)c^4+2m(1-n)c^2+m(m-1)}{c(1-c^2)^{3/2}}=0,
\end{aligned}\label{317}
\end{equation}
and
\begin{equation}
\begin{aligned}
&\frac{(n-m)[(1-n)c^2+m]}{c(1-c^2)^{3/2}}-\frac{n(n-1)(n-2)}{2}H_3\\
&=-\frac{(n(n-1)c^4+2m(1-n)c^2+m(m-1))(c^2(2n-1)-2m+1)}{c^3(1-c^2)^{3/2}}=0,
\end{aligned}\label{318}
\end{equation}
hence, from \eqref{00310}, \eqref{3-30}, \eqref{316}, \eqref{317}
and \eqref{318}, we have
\begin{equation}
\begin{aligned}
\lambda_1^{J_s}=-2n(n-1)H<\lambda_2^{J_s}=-n(n-1)H=\frac{n(n-1)(n-2)}{2}H_3<0.
\end{aligned}
\label{ex3}
\end{equation}
\end{ex}

In the following we will assume that $x:M\to \mathbb{S}^{n+1}(1)$ is
an n-dimensional compact orientable hypersurface with constant
scalar curvature $n(n-1)r,~r\geq 1$, in $\mathbb{S}^{n+1}(1),~n\geq
5$, when $r=1$, we assume moreover $H_3\neq 0.$ When $r>1$, we have
$n^2H^2>S>0$, when $r=1$, since $H_3\neq 0$, we have $H\neq 0$.
Hence, we can assume $H>0$ (cf. \cite{Cheng4} and \cite{HL}).

 Let $a$
be a fixed vector in $\mathbb{R}^{n+2}$. We define functions
$f^a:M\to \mathbb{R}$ and $\tilde{g}^a:M\to \mathbb{R}$ by
\begin{equation}
f^a=<a,x>,~\tilde{g}^a=<a,e_{n+1}>,
\end{equation}
where $x$ is the position vector and $e_{n+1}$ is the unit normal
vector.

By using the structure equations and the definition of the covariant
derivatives, we have the following result.

\begin{lem}(see \cite{CL})
The gradient and the second derivative of the functions $f$ and
$\tilde{g}$ are given by
\begin{equation}
\begin{aligned}
&f^a_i=<a,e_i>,~f^a_{ij}=\tilde{g}^ah_{ij}-f^a\delta_{ij},\\
&\tilde{g}^a_j=-\sum_{i=1}^{n}<a,e_i>h_{ij},~\tilde{g}^a_{jk}=-\sum_{i=1}^{n}<a,e_i>h_{ijk}-\sum_{i=1}^{n}\tilde{g}^ah_{ij}h_{ik}+f^ah_{jk}.
\end{aligned}\label{3-36}
\end{equation}
\end{lem}
\begin{proof}
By \eqref{2.1} we have
$$df^a=<a,dx>=\sum_{i}<a,e_i>\omega_i,$$
thus from \eqref{2.8} we have
\begin{equation}
f^a_i=<a,e_i>.\label{fi}
\end{equation}
From \eqref{2.2} and \eqref{fi} we have
$$
\begin{aligned}
\sum\limits_{j=1}^{n}f^a_{ij}\omega_j&=df_i+\sum\limits_{j=1}^{n}f_j\omega_{ji}=<a,de_i>+\sum_{j=1}^{n}<a,e_j>\omega_{ji}\\
&=\sum\limits_{j=1}^{n}<a,e_{n+1}>h_{ij}\omega_j-<a,x>\omega_i,
\end{aligned}
$$
hence we have
\begin{equation}
f^a_{ij}=<a,e_{n+1}>h_{ij}-<a,x>\delta_{ij}=\tilde{g}^ah_{ij}-f^a\delta_{ij}.\label{fij}
\end{equation}
After an analogous argument, we have
\begin{equation}
\tilde{g}^a_j=-\sum_{i=1}^{n}<a,e_i>h_{ij},~\tilde{g}^a_{jk}=-\sum_{i=1}^{n}<a,e_i>h_{ijk}-\sum_{i=1}^{n}\tilde{g}^ah_{ij}h_{ik}+f^ah_{jk}.
\end{equation}
\end{proof}

We will use a technique which was introduced by Li and Yau in
\cite{LY} and was later used by other authors (see \cite{MR},
\cite{P} and \cite{U}).

Let $B^{n+2}$ be the open unit ball in $\mathbb{R}^{n+2}$. For each
point $g\in B^{n+2}$, we consider the map
\begin{equation}
F_g(p)=\frac{p+(\mu<p,g>+\lambda)g}{\lambda(<p,g>+1)},~\forall~ p\in
\mathbb{S}^{n+1}(1)\subset\mathbb{R}^{n+2},
\end{equation}
where $\lambda=(1-\|g\|^2)^{-1/2}$, $\mu=(\lambda-1)\|g\|^{-2}$ and
$<,>$ denotes the usual inner product on $\mathbb{R}^{n+2}$. A
direct computation (see \cite{MR}, \cite{U}) shows that $F_g$ is a
conformal transformation from $\mathbb{S}^{n+1}(1)$ to
$\mathbb{S}^{n+1}(1)$ and the differential map $dF_g$ of $F_g$ is
given by
$$dF_g(v)=\lambda^{-2}(<p,g>+1)^{-2}\{\lambda(<p,g>+1)v-\lambda<v,g>p+<v,g>(1-\lambda)\|g\|^{-2}g\},$$
where $v$ is a tangent vector to $\mathbb{S}^{n+1}$ at the point
$p$. Hence, for two vectors $v,~w\in T_p\mathbb{S}^{n+1}$ we have
(see \cite{MR}, \cite{P} and \cite{U})
$$<dF_g(v),dF_g(w)>=\frac{1-\|g\|^2}{(<p,g>+1)^2}<v,w>.$$

By use of the technique in Li-Yau \cite{LY}, we have the following
result:
\begin{lem}(see \cite{MR}, \cite{P} and
\cite{U})

Let $x:M\to \mathbb{S}^{n+1}$ be a compact hypersurface in
$\mathbb{S}^{n+1}$ with constant scalar curvature $n(n-1)r,~r\geq1$,
and $u$ be a positive first eigenfunction of the Jacobi operator
$J_s$ on $M$, then there exists $g\in B^{n+2}$ such that
$\int_Mu(F_g\circ x)dv=(0,\ldots,0).$
\end{lem}

Let $\{E^A\}^{n+2}_{A=1}$ be a fixed orthonormal basis of
$\mathbb{R}^{n+2}$, for a fixed point $g\in B^{n+2}$, we define
functions $f^A:M\to \mathbb{R}( 1\leq A\leq n+2)$ by
\begin{equation}
f^A=<E^A,F_g\circ
x>=\frac{<E^A,x>+(\mu<x,g>+\lambda)<g,E^A>}{\lambda(<x,g>+1)},~\forall
1\leq A\leq n+2.\label{3.4}
\end{equation}

\begin{lem}
The gradient of $f^A$ is given by
\begin{equation}
f^A_i=\frac{<E^A,e_i>}{\lambda(<x,g>+1)}+\frac{<g,e_i>}{\lambda(<x,g>+1)^2}(-<E^A,x>+\frac{1-\lambda}{\lambda\|g\|^2}<g,E^A>).
\end{equation}
\end{lem}
\begin{proof}
By applying Lemma 3.3, we have
$$
\begin{aligned}
f^A_i&=\frac{<E^A,e_i>+\mu<g,e_i><g,E^A>}{\lambda(<x,g>+1)}-f^A\frac{<g,e_i>}{<x,g>+1}
\\
&=\frac{<E^A,e_i>}{\lambda(<x,g>+1)}+\frac{<g,e_i>}{\lambda(<x,g>+1)^2}(\mu<g,E^A>-<E^A,x>-\lambda<g,E^A>)\\
&=\frac{<E^A,e_i>}{\lambda(<x,g>+1)}+\frac{<g,e_i>}{\lambda(<x,g>+1)^2}(-<E^A,x>+\frac{1-\lambda}{\lambda\|g\|^2}<g,E^A>).
\end{aligned}
$$
\end{proof}

We also need the following Lemma 3.6, Lemma 3.7 and Lemma 3.8 to
estimate the second eigenvalue $\lambda_2^{J_s}$ of the Jacobi
operator $J_s$ on $M$.

\begin{lem}Let $M$ be an n-dimensional compact hypersurface with
constant scalar curvature $n(n-1)r,~r\geq 1$, in
$\mathbb{S}^{n+1}(1)$. Let $f^A$ be the function given by
\eqref{3.4}, we have
\begin{equation}
\sum_{A=1}^{n+2}\int_M(J_sf^A\cdot
f^A)dv=\int_M\frac{n(n-1)H(1-\|g\|^2)}{(<x,g>+1)^2}dv-\int_M\{\frac{n(n-1)}{2}(2H-(n-2)H_3+nHH_2)\}dv.\label{3.9}
\end{equation}
\end{lem}
\begin{proof}
By divergence theorem and Lemma 3.5 we have
\begin{equation}
\begin{aligned}
&-\sum_{A=1}^{n+2}\int_M(\Box f^A\cdot f^A)dv=\sum_{A=1}^{n+2}\int_M\sum_{i,j}(nH\delta_{ij}-h_{ij})f^A_if^A_jdv\\
&=\sum_{A=1}^{n+2}\int_M\sum_{i,j}(nH\delta_{ij}-h_{ij})(\frac{<E^A,e_i>}{\lambda(<x,g>+1)}+\frac{<g,e_i>}
{\lambda(<x,g>+1)^2}(-<E^A,x>+\frac{1-\lambda}{\lambda\|g\|^2}<g,E^A>))\\
&\cdot(\frac{<E^A,e_j>}{\lambda(<x,g>+1)}+\frac{<g,e_j>}{\lambda(<x,g>+1)^2}(-<E^A,x>+\frac{1-\lambda}{\lambda\|g\|^2}<g,E^A>))dv\\
&=\int_M\{\sum_{i,j}[nH\delta_{ij}-h_{ij}][\frac{\delta_{ij}}{\lambda^2(<x,g>+1)^2}+\frac{<g,e_i><g,e_j>}{\lambda^4\|g\|^2(<x,g>+1)^2}[
2(1-\lambda)\lambda(<x,g>+1)\\
&+\lambda^2\|g\|^2-2(1-\lambda)\lambda<x,g>+(1-\lambda)^2]
]\}dv\\
&=\int_M\sum_{i,j}(nH\delta_{ij}-h_{ij})\cdot\frac{\delta_{ij}}{\lambda^2(<x,g>+1)^2}dv\\
&=\int_M\frac{n(n-1)H(1-\|g\|^2)}{(<x,g>+1)^2}dv,
\end{aligned}
\label{3.10}
\end{equation}
where we use the fact that
$\sum\limits_{A=1}^{n+2}<E_A,X><E_A,Y>=<X,Y> (\forall~ X, Y \in
\mathbb{R}^{n+2})$ in the third equality.

By Newton formula, we have
\begin{equation}
\begin{aligned}
f_3&=n^3H^3+\frac{n(n-1)(n-2)}{2}H_3-\frac{3n^2(n-1)}{2}HH_2,\\
S&=n^2H^2-n(n-1)H_2.
\end{aligned}\label{3-45}
\end{equation}
Thus $J_s$ becomes
\begin{equation}
\begin{aligned}
J_s&=-\Box-\{n(n-1)H+nH(n^2H^2-n(n-1)H_2)\\
&-(n^3H^3+\frac{n(n-1)(n-2)}{2}H_3-\frac{3n^2(n-1)}{2}HH_2)\}\\
&=-\Box-n(n-1)H-\frac{n^2(n-1)}{2}HH_2+\frac{n(n-1)(n-2)}{2}H_3\\
&=-\Box-\frac{n(n-1)}{2}(2H-(n-2)H_3+nHH_2).
\end{aligned}
\label{1.4}
\end{equation}
Then by using the fact that
\begin{equation}\sum_{A=1}^{n
+2}f^A\cdot f^A=\sum_{A=1}^{n+2}<E^A,F_g\circ x><E^A,F_g\circ
x>=<F_g\circ x,F_g\circ x>=1,\label{319}
\end{equation} we  immediately get \eqref{3.9}.
\end{proof}

For a fixed point $g\in B^{n+2}$, let
\begin{equation}f=<x,g>,~\tilde{g}=<e_{n+1},g>,~\rho=-\ln{\lambda}-\ln{(1+f)},\label{3.19}\end{equation}
where $\lambda=(1-\|g\|^2)^{-1/2}$, $x$ is the position vector and
$e_{n+1}$ is the unit normal vector. We have
\begin{equation}e^{2\rho}=\frac{1}{\lambda^2(1+f)^2}=\frac{1-\|g\|^2}{(<x,g>+1)^2},~\rho_i=\frac{-f_i}{1+f},~\rho_{ij}
=\frac{-f_{ij}}{1+f}+\frac{f_if_j}{(1+f)^2}.\label{03-34}\end{equation}

\begin{lem}Let $x:M\to \mathbb{S}^{n+1}(1)$ be an n-dimensional compact hypersurface with constant scalar curvature
$n(n-1)r,~r\geq 1$, in $\mathbb{S}^{n+1}(1)$. When $r=1$, we assume
moreover that $H_3\neq 0.$ Then we have $H\neq 0$, hence we can
assume $H>0$. Let $\rho$ be the function defined by \eqref{3.19}, we
have
\begin{equation}
\int_M\frac{H(1-\|g\|^2)}{(<x,g>+1)^2}dv\leq\int_M(H+\frac{H_2^2}{H})dv-\int_M[H\|\nabla
\rho\|^2-\frac{2}{n(n-1)}\sum_{i,j}(nH\delta_{ij}-h_{ij})\rho_i\rho_j]dv,\label{3.11}
\end{equation}
and the equality holds if and only if
$H_2+\frac{\tilde{g}H}{1+f}\equiv0$ on $M$.
\end{lem}

\begin{proof} Under the hypothesis of the lemma, we can assume $H>0$ (cf. \cite{
Cheng4} and \cite{ABS}). We have
\begin{equation}
\sum_{i,j}(nH\delta_{ij}-h_{ij})\rho_i\rho_j=\sum_{i,j}(nH\delta_{ij}-h_{ij})\frac{f_if_j}{(1+f)^2}=\frac{nH\|\nabla
f\|^2}{(1+f)^2}-\sum_{i,j}\frac{h_{ij}f_if_j}{(1+f)^2}, \label{3.12}
\end{equation}
and
\begin{equation}
\begin{aligned}
\Box
\rho&=\sum_{i,j}(nH\delta_{ij}-h_{ij})\rho_{ij}=\sum_{i,j}(nH\delta_{ij}-h_{ij})(\frac{-f_{ij}}{1+f}+\frac{f_if_j}{(1+f)^2})\\
&=\frac{-\Delta fnH}{1+f}+\frac{nH\|\nabla f\|^2}{(1+f)^2}
+\sum_{i,j}\frac{h_{ij}f_{ij}}{1+f}-\sum_{i,j}\frac{h_{ij}f_if_j}{(1+f)^2}.
\label{3.13}
\end{aligned}
\end{equation}

From \eqref{03-34}, \eqref{3.12} and \eqref{3.13} and  by using
Lemma 3.3, we have
$$
\begin{aligned}
&(\Box
\rho-\sum_{i,j}(nH\delta_{ij}-h_{ij})\rho_i\rho_j)\cdot\frac{2}{n(n-1)}+\frac{H(1-\|g\|^2)}{(1+f)^2}\\\
&=(\frac{-\Delta
fnH}{1+f}+\sum_{i,j}\frac{h_{ij}f_{ij}}{1+f})\cdot\frac{2}{n(n-1)}+\frac{H(1-\|g\|^2)}{(1+f)^2}\\
&=(\frac{-nH(nH\tilde{g}-nf)}{1+f}+\sum_{i,j}\frac{h_{ij}(\tilde{g}h_{ij}-f\delta_{ij})}{1+f})\cdot\frac{2}{n(n-1)}+
\frac{H(1-\|g\|^2)}{(1+f)^2}\\
&=\frac{2Hf-2H_2\tilde{g}}{1+f}+\frac{H(1-f^2-\sum\limits_{i}f_i^2-\tilde{g}^2)}{(1+f)^2}=H-\sum\limits_{i}\frac{Hf_i^2}{(1+f)^2}-\frac{H\tilde{g}^2}{(1+f)^2}-\frac{
2H_2\tilde{g}}{1+f}\\
&=H-\sum\limits_{i}\frac{Hf_i^2}{(1+f)^2}+\frac{H_2^2}{H}-\frac{(H_2+\frac{\tilde{g}H}{1+f})^2}{H}=H+\frac{H_2^2}{H}-H\|\nabla
\rho\|^2-\frac{(H_2+\frac{\tilde{g}H}{1+f})^2}{H},
\end{aligned}
$$
which immediately implies
\begin{equation}
\begin{aligned}
&\int_M\frac{H(1-\|g\|^2)}{(<x,g>+1)^2}dv\\
&=\int_M[H+\frac{H_2^2}{H}-H\|\nabla
\rho\|^2+\frac{2}{n(n-1)}\sum\limits_{i,j}(nH\delta_{ij}-h_{ij})\rho_i\rho_j-\frac{(H_2+\frac{\tilde{g}H}{1+f})^2}{H}]dv.
\end{aligned}
\end{equation}
Hence we get the inequality \eqref{3.11} and the equality holds if
and only if $H_2+\frac{\tilde{g}H}{1+f}\equiv0$ on $M$.\end{proof}

\begin{lem}Let $M$ be an n-dimensional compact hypersurface with
constant scalar curvature $n(n-1)r,~r\geq 1$, in
$\mathbb{S}^{n+1}(1),~n\geq5$. When $r=1$, we assume moreover that
$H_3\neq 0.$ Then we have $H\neq 0$, hence we can assume $H>0$. We
have
\begin{equation}\int_M[H\|\nabla
\rho\|^2-\frac{2}{n(n-1)}\sum\limits_{i,j}(nH\delta_{ij}-h_{ij})\rho_i\rho_j]dv\geq
0.\label{3.24}\end{equation}
\end{lem}

\begin{proof}Under the hypothesis of the lemma, we can assume $H>0$ (cf. \cite{
Cheng4} and \cite{ABS}). $\forall ~p\in M$, let $k_1,\ldots,k_n$
denote the principal curvatures of $M$ at $p$, we choose an
orthonormal basis such that $h_{ij}=\delta_{ij}k_i$. By Gauss
equation \eqref{2.6}, we have
\begin{equation}n^2H^2-\sum_{i}k_i^2=n(n-1)(r-1)\geq 0,\label{3.25}\end{equation}
which leads to
\begin{equation}nH\geq |k_i|,~\forall 1\leq i\leq n.\label{3.26}\end{equation}

As $n\geq 5$, we have $\frac{n(n-3)}{2}H\geq nH$, so we have
$$
\begin{aligned}&H\|\nabla
\rho\|^2-\frac{2}{n(n-1)}\sum\limits_{i,j}(nH\delta_{ij}-h_{ij})\rho_i\rho_j\\&=H\sum\limits_{i}
\rho_i^2-\frac{2}{n(n-1)}\sum\limits_{i,j}(nH\delta_{ij}-\delta_{ij}k_i)\rho_i\rho_j\\
&=H\sum\limits_{i} \rho_i^2-\sum\limits_{i}\frac{2}{n(n-1)}(nH-k_i)\rho_i^2\\
&=\frac{2}{n(n-1)}\sum\limits_{i}
\rho_i^2(\frac{n(n-3)}{2}H+k_i)\geq\frac{2}{n(n-1)} \sum\limits_{i}
\rho_i^2( nH-|k_i|)\geq 0.
\end{aligned}$$
Hence, we get $H\|\nabla
\rho\|^2-\frac{2}{n(n-1)}\sum\limits_{i,j}(nH\delta_{ij}-h_{ij})\rho_i\rho_j\geq
0$ holds at every point of $M$, which immediately implies
\eqref{3.24}.
\end{proof}

\section{Proofs of Theorem 1.3 and Theorem 1.4}
\noindent{\bf Proof of Theorem 1.3:} Since $r>1$, we have $\Box$ is
an elliptic operator and $H\neq 0$. Hence, we can assume $H>0$ (see
\cite{Cheng4}). Let $u$ be a first eigenfunction of $J_s$, we can
assume $u$ is positive on $M$, by Lemma 3.4 there exists $g\in
B^{n+2}$ such that
\begin{equation}
\int_Mu (F_g\circ x)dv=(0,\ldots,0),\label{4.1}
\end{equation}
which implies that the functions $\{f^A,~1\leq A\leq n+2\}$ given by
\eqref{3.4} are perpendicular to the function $u$, i.e.,
$\int_Mu\cdot f^Adv=0,~\forall 1\leq A\leq n+2$. Then by using the
min-max characterization of eigenvalues for elliptic operators, we
have
\begin{equation}
\lambda_2^{J_s}\cdot\int_M(f^A\cdot f^A)dv\leq \int_M(J_sf^A\cdot
f^A)dv,~\forall~ 1\leq A\leq n+2.\label{4.2}
\end{equation}
Summing up and using the fact that $\sum\limits_{A=1}^{n+2}f^A\cdot f^A=1$
(see \eqref{319}), we obtain
\begin{equation}
\lambda_2^{J_s}\cdot Vol(M) \leq \sum_{A=1}^{n+2}\int_M(J_sf^A\cdot f^A)dv.
\label{3.16}
\end{equation}
From Lemma 3.6 and \eqref{3.16} we have
\begin{equation}
\lambda_2^{J_s}\cdot Vol(M)\leq\
\int_M\frac{n(n-1)H(1-\|g\|^2)}{(<x,g>+1)^2}dv-\int_M\frac{n(n-1)}{2}(2H-(n-2)H_3+nHH_2)dv.\label{3.17}
\end{equation}

Then by \eqref{3.17}, Lemma 3.7 and Lemma 3.8, we have
\begin{equation}
\begin{aligned}
\lambda_2^{J_s}\cdot Vol(M) &\leq
n(n-1)\cdot\int_M(H+\frac{H_2^2}{H})dv-\int_M\frac{n(n-1)}{2}(2H-(n-2)H_3+nHH_2)dv\\
&=
n(n-1)\cdot\int_M(\frac{H_2^2}{H}+\frac{n-2}{2}H_3-\frac{nHH_2}{2})dv.
\end{aligned}\label{4.5}
\end{equation}

From definition of $H_2$ and the Gauss equation \eqref{2.6} we have
\begin{equation}H_2=r-1=\text{constant}>0.\label{4-6}\end{equation} So we have $H_3\leq
\frac{H_2^2}{H}$ and $H_2\leq H^2$ (see \cite{HLP}, p. 52) and hence
\begin{equation}
\begin{aligned}
\lambda_2^{J_s} \cdot Vol(M)&\leq
n(n-1)\cdot\int_M(\frac{H_2^2}{H}+\frac{n-2}{2}H_3-\frac{nHH_2}{2})dv\\
&\leq
n(n-1)\cdot\int_M(\frac{H_2^2}{H}+\frac{n-2}{2}\frac{H_2^2}{H}-\frac{nHH_2}{2})dv\\
&= n(n-1)\cdot\int_M\frac{nH_2}{2}(\frac{H_2}{H}-H)dv\leq 0,
\end{aligned}
\label{4.6}
\end{equation}
therefore we get $\lambda_2^{J_s}\leq 0$.

When $\lambda_2^{J_s}=0$, then all the inequalities become
equalities. From \eqref{4.6} we have $H_2=H^2$ on $M$, since $H_2$
is a positive constant, we get $M$ is a totally umbilical and
non-totally geodesic hypersurface with constant scalar curvature
$n(n-1)r$. On the other hand, if $M$ is a totally umbilical and non-totally geodesic
hypersurface with constant scalar curvature $n(n-1)r$, from Example
3.1 in section 3, we know that $\lambda_2^{J_s}=0.$\qed

\begin{rem}
We notice that from \eqref{4.6} we can get a more precise upper
bound for $\lambda_2^{J_s}$, that is,
\begin{equation}
\begin{aligned}
\lambda_2^{J_s}&\leq
n(n-1)(\frac{H_2^2}{\min{H}}+\frac{n-2}{2}\max{H_3}-\frac{nH_2}{2}\min{H})\\
&=n(n-1)(\frac{(r-1)^2}{\min{H}}+\frac{n-2}{2}\max{H_3}-\frac{n(r-1)}{2}\min{H})
.
\end{aligned}
\end{equation}
\end{rem}

\noindent{\bf Proof of Theorem 1.4:} Since $r=1$, from \eqref{4-6}
we have $H_2=0$. Since we assume that $H_3$ does not vanish on $M$,
we have $J_s$ is elliptic and the mean curvature $H$ does not vanish
on $M$(cf. Proposition 1.5 in \cite{HL}). Hence, we can assume
$H>0$. Thus $H_3\leq \frac{H_2^2}{H}=0$. Since we assume that
$H_3\neq0$ on $M$, we get $H_3<0.$ As Lemma 3.6, Lemma 3.7 and Lemma
3.8 hold for both the case $r>1$ and the case $r=1$, after an
analogous argument with the proof of Theorem 1.3, we know that
\eqref{4.1}-\eqref{4.5} still hold in this case, hence we have
\begin{equation}
\begin{aligned}
\lambda_2^{J_s}\cdot Vol(M) &\leq
n(n-1)\cdot\int_M(\frac{H_2^2}{H}+\frac{n-2}{2}H_3-\frac{nHH_2}{2})dv\\
&=\frac{
n(n-1)(n-2)}{2}\cdot\int_MH_3dv\\
&\leq \frac{ n(n-1)(n-2)}{2}\max{H_3}\cdot Vol(M)\\
&=-\frac{ n(n-1)(n-2)}{2}\min{|H_3|}\cdot Vol(M).
\end{aligned}\label{4.7}
\end{equation}
Hence, we get \begin{equation}\lambda_2^{J_s}\leq-\frac{
n(n-1)(n-2)}{2}\min{|H_3|}.\label{4.9}\end{equation}

When $\lambda_2^{J_s}=-\frac{ n(n-1)(n-2)}{2}\min{|H_3|}$,  the
inequalities in \eqref{3.11}, \eqref{4.2} and \eqref{4.7} become
equalities. The equality in \eqref{4.7} holds implies that
$H_3=\text{constant}\neq 0$. Since $H_2=0$, the equalities in
\eqref{3.11} holds implies that $\tilde{g}=<g,e_{n+1}>\equiv 0$ on
$M$. We claim that $g$ must be $0$, otherwise, we have that $M$ is a
hypersphere (see Theorem 1 in \cite{NS}), hence $M$ is totally
umbilical, since $H_2=0$, we immediately get $M$ is totally geodesic
which is a contradiction with $H_3\neq 0$. Hence we have $g\equiv0$,
from \eqref{3.4} we get $f^A=<E^A,F_g\circ x>=<E^A,x>$, which means
$\{f^A,~1\leq A\leq n+2\}$ are the position functions of
 $x:M\to\mathbb{S}^{n+1}(1)$.
Since the equality in \eqref{4.2} holds, it follows  that the
position functions $\{f^A=<E^A,x>,~1\leq A\leq n+2\}$ must be the
second eigenfunctions of $J_s$ corresponding to $\lambda_2^{J_s}$.

On the other hand, if we assume that $H_3=\text{constant}\neq 0$ and
the position functions $\{\tilde{f}^A=<E^A,x>,~1\leq A\leq n+2\}$
are the second eigenfunctions of $J_s$ corresponding to
$\lambda_2^{J_s}$. Since  $H_3\neq 0$, we have $H\neq 0$. Hence, we
can assume $H>0$, $H_3<0$ (cf. Proposition 1.5 in \cite{HL}).

 Since
$H_2=0$, by using \eqref{1-1} and \eqref{3-36}, we get $$\Box
\tilde{f}^A=n(n-1)H_2<E^A,e_{n+1}>-n(n-1)H\tilde{f}^A=-n(n-1)H\tilde{f}^A,~\forall
~1\leq A\leq n+2,$$ then from \eqref{1.1} and \eqref{3-45} we have
$$
\begin{aligned}
J_s\tilde{f}^A&=n(n-1)H\tilde{f}^A-\{n(n-1)H+nHS-f_3\}\tilde{f}^A\\
&=(f_3-nHS)\tilde{f}^A\\
&=\{(n^3H^3+\frac{n(n-1)(n-2)}{2}H_3-\frac{3n^2(n-1)}{2}HH_2)-(n^3H^3-n^2(n-1)HH_2)\}\tilde{f}^A\\
&=\frac{n(n-1)(n-2)}{2}H_3\tilde{f}^A,  ~\forall ~1\leq A\leq n+2,
\end{aligned}
$$
hence we get $\lambda_2^{J_s}=\frac{n(n-1)(n-2)}{2}H_3=-\frac{
n(n-1)(n-2)}{2}\min{|H_3|}$.

In particular, when $M$ is the Riemannian product
$\mathbb{S}^{m}(c)\times\mathbb{S}^{n-m}(\sqrt{1-c^2}),~1\leq m\leq
n-2$ with $c=\sqrt{\frac{(n-1)m+\sqrt{(n-1)m(n-m)}}{n(n-1)}}$, from
Example 3.3 in section 3, we know that
 the equality in \eqref{4.9} is attained.
 \qed

\begin{rem}
Since Lemma 3.8 does not hold when $n=3$ and $n=4$, we can not prove
Theorem 1.3 and Theorem 1.4 by our technique in $n=3$ and $n=4$. So
it is an interesting problem to study the estimate for the second
eigenvalue of the Jacobi oeprator $J_s$ of the hypersurface
$x:M^n\to \mathbb{S}^{n+1}(1)$ when $n=3$ and $n=4$.
\end{rem}

\vskip 1cm
\begin{flushleft}
\medskip\noindent

Haizhong Li: {\sc Department of Mathematical Sciences, Tsinghua
University, Beijing 100084, People's Republic of China} \ \ E-mail:
hli@math.tsinghua.edu.cn

Xianfeng Wang: {\sc Department of Mathematical Sciences, Tsinghua
University, Beijing 100084, People's Republic of China} \ \ E-mail:
xf-wang06@mails.tsinghua.edu.cn

\end{flushleft}
% ------------------------------------------------------------------------
\end{document}